\documentclass[11pt,twoside]{article}
\usepackage{amsmath,amssymb,latexsym}

\newtheorem{theorem}{Theorem}[section]
\newtheorem{lemma}[theorem]{Lemma}
\newtheorem{corollary}[theorem]{Corollary}

\newtheorem{definition}[theorem]{Definition}
\newtheorem{remark}[theorem]{Remark}
\newtheorem{proposition}[theorem]{Proposition}
\newtheorem{claim}[theorem]{Claim}

\newtheorem{construction}[theorem]{Construction}

\newcommand{\pproof}{\quad{\sc Proof of  }}

\newcommand{\proof}{\quad{\sc Proof}:  \,}
\newcommand{\bewend}{\vspace{1ex}\hfill $ \Box $}

\newcommand{\OO}{{\cal O}}
\newcommand{\LL}{{\cal L}}
\newcommand{\BB}{{\cal B}}
\newcommand{\EE}{{\cal E}}

\newcommand{\MM}{{\cal M}}
\newcommand{\GG}{{\cal G}}

\newcommand{\KK}{{\cal K}}
\newcommand{\CC}{{\cal C}}
\newcommand{\HH}{{\cal H}}
\newcommand{\II}{{\cal I}}

\renewcommand{\P}{\mathbb P}
\newcommand{\N}{\mathbb N}
\newcommand{\Q}{\mathbb Q}
\newcommand{\C}{\mathbb C}
\newcommand{\JJ}{{\cal J}}
\newcommand{\Z}{\mathbb Z}

\newcommand{\oideal}[3]{\omega_{#1}\left\{-\frac{#2}{#3}\right\}}
\newcommand{\noideal}[3]{\omega_{#1}\left(-\left[\frac{#2}{#3}\right]\right)}
\newcommand{\noidealtext}[3]{\omega_{#1}\left(-\left[\frac{#2}{#3}\right]\right)}

\newcommand{\oidealtext}[3]{\omega_{#1}\{-\frac{#2}{#3}\}}

\newcommand{\coker}{{\rm coker}}

\newcommand{\codim}{{\rm codim}}

\begin{document}

\thispagestyle{empty}

\begin{center}
{\large \sc WEAK POSITIVITY AND DYSON'S LEMMA} 

{\small \sc Markus Wessler}
\end{center}

\begin{abstract}
We give a proof of Dyson's Lemma for a product of smooth projective
varieties of arbitrary dimension.
\end{abstract}

\


\begin{center}
\section{\sc Introduction}
\end{center}

By Liouville's theorem, complex numbers which can be approximated
by rational numbers very well are necessarily transcendent. 
In other words, 
if $\alpha$ is algebraic of degree $d \geq 2$, then
for all $\varepsilon > 0$ there are only finitely many rational
numbers $p/q$ such that
$|p/q - \alpha| \leq q^{-(d + \varepsilon)}$.
But since in the situation above an infinite sequence of
rational numbers $(p_n/q_n)_{n \in \N}$ could be constructed, 
satisfying $|p_n/q_n - \alpha| \leq q_n^{-2}$ for all
$n \in \N$, 
it was clear that there was some lower bound for the exponent. 
It took about hundred years until Roth in 1955 could 
prove that replacing $d$ by $2$ was in fact the optimal bound.

In most approaches, auxiliary polynomials in two or more variables
were constructed. In his paper \cite{d47}, Dyson describes explicitly 
which properties these polynomials should have. 
In order to replace the exponent $d + \varepsilon$
by $\sqrt{2d} + \varepsilon$, he proved a statement about the
existence of certain polynomials, which is
known as Dyson's Lemma today.

Fixing a certain number of points in the complex plain, Dyson
considers polynomials in two variables of multidegree $\underline{d}
= (d_1, d_2)$ with respect to a special kind of zero conditions 
in these points, called the index.
Clearly there exists some polynomial satisfying them, provided the 
number of conditions
is bounded by $d_1 \cdot d_2$. Now Dyson finds out that,
increasing the number of conditions too much, it becomes
impossible to find such a polynomial at all. 
In other words, assuming the
existence of such a polynomial, he shows that the number of conditions
is necessarily bounded by $C \cdot d_1 \cdot d_2$ for some constant
$C$ depending on $d_1$ and $d_2$.
Moreover, this constant tends to one when increasing the ratio
$d_2/d_1$, which means that the number of conditions is asymptotically 
independent.

We want to reformulate this in terms of algebraic geometry. 
The question should be the following. If one compactifies the 
situation by considering $(\P^1)^n$ over $\C$ (or over some 
algebraically closed field) and
if one identifies polynomials with sections of
some special sheaf and encodes the zero conditions into some ideal sheaf,
how can the existence of a special section be interpreted then? This
immediately leads to the notion of weak positivity, and it is exactly
what Esnault-Viehweg stated and proved
in \cite{ev84}, using positivity statements and vanishing theorems.

The next question arising is: how can this be generalized further?
The case of a product of arbitrary curves has been treated by
Nakamaye \cite{n95}, for example, using derivations. 
The proof of Dyson's Lemma in \cite{ev84}, however, is
based on positivity methods from algebraic geometry.
Trying to give a corresponding proof for a product of arbitrary 
curves, one is led to a much more general situation. We shall, 
therefore, in section 5 give a proof for
Dyson's Lemma as it is stated below (\ref{dysonv}) and then, in
section 6, deduce from this the result for a product of curves, where
a slightly stronger positivity statement can be obtained.

We now consider the following situation.
Let $X = X_1 \times \dots \times X_n$ be a product of smooth 
projective 
varieties, defined over an algebraically closed field $k$ of
characteristic zero. 
For $\nu = 1 , \dots , n$, let us denote by $m_{\nu}$ the dimension of
$X_{\nu}$ and by $p_{\nu} : X \to X_{\nu}$ the projection onto the
$\nu$-th factor.
Let us fix very ample sheaves $\LL_{\nu}$ on $X_{\nu}$ and let us
write $e_{\nu} = c_1(\LL_{\nu})^{m_{\nu}}$.
For every $n$-tuple $\underline{\delta} = (\delta_1 , \dots,
\delta_n)$ of non-negative integers
let us
denote by $\underline{\LL}^{\underline{\delta}} = \LL_1^{\delta_1}
\boxtimes \dots \boxtimes \LL_n^{\delta_n}$ the induced sheaf on $X$.

Let $S = \{\xi_1 , \dots , \xi_M\}$ be a finite set of
$M \geq 2$ points on $X$ with projections 
$\xi_{\mu, \nu} = p_{\nu}(\xi_{\mu})$ onto the factors.
Let us fix positive rational
numbers $t_1, \dots , t_M$ and an $n$-tuple $\underline{d} =
(d_1 , \dots , d_n)$ of non-negative integers. 
Let us assume that $d_1 \cdot e_1 \geq \dots \geq d_n \cdot e_n$.
For $\mu = 1 , \dots , M$ we define 
$\II_{\xi_{\mu}}$ to be the ideal sheaf generated by 
$ m_{\xi_{\mu,1}}^{\alpha_1} \boxtimes \dots
\boxtimes m_{\xi_{\mu,n}}^{\alpha_n}$ for $\frac{\alpha_1}{d_1} +
\dots + \frac{\alpha_n}{d_n} \geq t_{\mu}$ and 
$\II = \bigcap_{\mu = 1}^M \II_{\xi_{\mu} , t_{\mu}}.$

Finally, let us fix non-negative integers $\gamma_{\nu}$ such 
that $\LL_{\nu}^{\gamma_{\nu}} \otimes
\omega_{X_{\nu}}^{-1}$ is globally generated and let us define
$M'_{\nu} = \min \{2, |p_{\nu}(S)|\}$.
In this situation, we are going to prove:

\begin{theorem}\label{dysonv} 
If $\underline{\LL}^{\underline{d}} \otimes \II$ is
effective, then $\underline{\LL}^{\underline{d'}} \otimes \II$ is
weakly positive over a product open set, where $\underline{d'} =
(d_1' , \dots , d_n')$ with $$d_{\nu}' = d_{\nu} + M_{\nu} \cdot
\sum_{j=\nu+1}^n d_j \cdot e_j$$ and $M_{\nu} = M_{\nu}' \cdot
m_{\nu} + \gamma_{\nu}$ for $\nu = 1 , \dots , n$.
\end{theorem}

Comparing this to the classical situation, we notice that the ideal sheaf
$\II$ carries the index conditions, the effectivity of 
$\underline{\LL}^{\underline{d}} \otimes \II$ corresponds to the
existence of a polynomial of multidegree $\underline{d}$ satisfying
these conditions, and the maximal number of possible conditions is
encoded in $\underline{d'}$.


\begin{center}
\section{\sc Weak Positivity}
\end{center}

Let us recall the following positivity notions:

\begin{definition}
\rm
Let $X$ be a quasi-projective variety and let $U \subseteq X$ be an open 
subset.

\begin{enumerate}

\item A locally free sheaf $\GG$ on $X$ is called {\it weakly
positive over $U$}, if there exists an ample invertible sheaf $\HH$ on $X$
such that for all $\alpha > 0$ the sheaf $S^{\alpha}(\GG)
\otimes \HH$ is semi-ample over $U$, which means that for some $\beta
> 0$ the sheaf $S^{\alpha \beta}(\GG)
\otimes \HH^{\beta}$  is globally generated over $U$.
If $\GG$ is weakly positive over $X$, we call $\GG$ weakly positive.
Obviously, we may replace $\beta$ by arbitrary multiples.

\item An invertible sheaf $\LL$ on $X$ is called 
{\it very ample with respect to $U$}, if $\LL$ is globally
generated over $U$ by sections of a finite dimensional subspace $V
\subseteq H^0(X,\LL)$ and the natural map $U \to \P(V)$ defined by
these sections is an embedding.
$\LL$ is called {\it ample with respect to  $U$}, if some power 
of $\LL$ is very ample with respect to $U$.
\end{enumerate}
\end{definition}

\begin{lemma}\label{lokampkrit}
Let $X$ be a quasi-projective variety, let $\LL$ be an invertible sheaf
on $X$ and let
$U \subseteq X$ be an open subset. Then $\LL$ is ample with respect 
to $U$ if and only if there exists
a blowing up $\tau : X' \to X$ with centre outside $U$ such that
for some ample invertible sheaf $\LL'$ on $X'$ and some $\mu > 0$ 
we have an inclusion
$\LL' \to \tau^{\ast}\LL^{\mu}$, which is an isomorphism over
$\tau^{-1}(U)$.
\end{lemma}

\proof We only have to show the condition is necessary and to this end
we may assume that $\LL$ is very ample with respect to $U$. 
Let $V$ be the space of sections generating $\LL$ over $U$. Thus we
have a rational map $\varphi : X \to \P(V)$,
given by the map $V \otimes \OO_X \to \LL$, surjective over $U$.
If $\GG$ is the image sheaf of this map, then we consider the
blowing up $\tau' : X'' \to X$ with respect to the ideal sheaf 
$\GG \otimes \LL^{-1}$. By \cite{h77}, II.7.17.3 there exists a morphism
$\varphi'' : X'' \to \P(V)$ and an inclusion
${\varphi''}^{\ast}\OO_{\P(V)}(1) \to {\tau'}^{\ast}\LL$. Since
${\varphi''}^{\ast}\OO_{\P(V)}(1)$ is not necessarily ample, we
have to continue blowing up. Let $Y$ be the image of $X''$ under
$\varphi''$ and $M \in \N-\{0\}$ such that $X'' \subseteq \P^M$.
Then $\varphi''$ factorizes over $Y \times \P^M = \P(\EE)$, where
$\P(\EE)$ is the projective bundle of $\EE =
p_2^{\ast}\left(\bigoplus^{M+1}\OO_{\P^M}\right)$ (see \cite{h77},
II.7). We have natural maps $$\EE = {p_1}_{\ast}\OO_{\P(\EE)}(1)
\to \varphi''_{\ast}\left(\OO_{\P(\EE)}(1)|_ {X''}\right) \to
\OO_Y(1).$$ 
Let $\BB$ be the image sheaf of this composed map.
We consider the blowing up $\varphi' :
X' \to Y$ of $\BB \otimes \OO_Y(-1)$, which is an ideal sheaf, since 
$\varphi''$ is birational. The sheaf $$\BB' =
{\varphi'}^{-1}\left(\BB \otimes \OO_Y(-1)\right) \cdot \OO_{X'}$$ is
therefore invertible on $X'$, and we obtain $\BB' \otimes
\OO_{X'}(1)$ as an invertible quotient of ${\tau'}^{\ast}\EE$. By
\cite{h77}, II.7.12 this corresponds to a morphism $X' \to
\P(\EE)$, factorizing canonically over some morphism $\eta : X' \to
X''$. Then $\tau = \tau' \circ \eta$ is the blowing up we need.

By construction we obtain an exceptional divisor $E$ for $\tau$ such
that $\OO_{X'}(-E)$ is relatively ample. Then there exists some 
$\mu > 0$ such that
$$\LL' = \eta^{\ast}{\varphi''}^{\ast}\OO_Y(\mu) \otimes \OO_{X'}(-E)$$
is ample. This sheaf is contained in $\eta^{\ast}{\tau'}^{\ast}\LL^{\mu} =
\tau^{\ast}\LL^{\mu}$,
and this inclusion is an isomorphism over $\tau^{-1}(U)$.
\bewend

This immediately implies the compatibility of locally ample sheaves
with finite morphisms:

\begin{corollary}\label{lokampend}
Let $\sigma : Y \to X$ be a morphism of normal quasi-projective
varieties, let $U \subseteq X$
be an open subset such that $\sigma|_{\sigma^{-1}(U)}$
is finite, and let $\LL$ be an invertible sheaf on $X$.
If $\LL$ is ample with respect to $U$, then
$\sigma^{\ast}\LL$ is ample with respect to $\sigma^{-1}(U)$.
\end{corollary}

\proof If $\sigma$ is finite, then by \ref{lokampkrit}
we find a blowing up $\tau : X' \to X$ with centre
outside $U$, an ample sheaf $\LL'$ on $X'$, some $\mu > 0$ and an
inclusion $\LL'\to \tau^{\ast}\LL^{\mu}$, being an isomorphism
over $\tau^{-1}(U)$. Let
\begin{eqnarray*}
Y' & \stackrel{\sigma'}{\longrightarrow} & X' \\
\tau' \downarrow & & \downarrow \tau \\
Y & \stackrel{\sigma}{\longrightarrow} & X
\end{eqnarray*}
be the fibre product. Since $\sigma'$ is finite, ${\sigma'}^{\ast}\LL'$
is ample on $Y'$. Moreover, we have an inclusion
$${\sigma'}^{\ast}\LL' \longrightarrow
{\sigma'}^{\ast}\tau^{\ast}\LL^{\mu}
= {\tau'}^{\ast}\sigma^{\ast}\LL^{\mu} ,$$
being an isomorphism over $\sigma^{-1}(U)$, which by
\ref{lokampkrit} implies that $\sigma^{\ast}\LL$ is ample with respect to
$\sigma^{-1}(U)$.

For the general case we may, considering the Stein factorization
of $\sigma$,
assume that $\sigma$ is birational and an isomorphism over $U$. But
then we are done, since $\sigma_{\ast} \sigma^{\ast} \LL = \LL$ and
hence the sections of $\LL$ correspond to the sections of 
$\sigma^{\ast}\LL$.
\bewend

The definition of weak positivity is independent of
the choice of the ample invertible sheaf $\HH$. Moreover, we have:

\begin{theorem}\label{amplunabh}
Let $X$ a quasi-projective  variety, $\GG$ a locally free sheaf on
$X$ and $U \subseteq X$ an open subset. Then the following
statements are equivalent:

(a) $\GG$ is weakly positive over $U$.

(b) There exists an ample invertible sheaf $\HH$ on $X$ such that for
all $\alpha > 0$ the sheaf $S^{\alpha \beta}(\GG) \otimes \HH^{\beta}$
is globally generated over $U$ for some $\beta > 0$.

(c) For every ample invertible sheaf $\HH$ on $X$ and for
all $\alpha > 0$, the sheaf $S^{\alpha \beta}(\GG) \otimes \HH^{\beta}$
is globally generated over $U$ for some $\beta > 0$.

(d) There exists an invertible sheaf $\HH$ on $X$, ample with respect to
$U$, such
that for all $\alpha > 0$ the sheaf $S^{\alpha \beta}(\GG) \otimes
\HH^{\beta}$
is globally generated over $U$ for some $\beta > 0$.

(e) For every invertible sheaf $\HH$ on $X$, ample with respect to $U$,
and for
all $\alpha > 0$, the sheaf $S^{\alpha \beta}(\GG) \otimes \HH^{\beta}$
is globally generated over $U$ for some $\beta > 0$.

(f) There exists an invertible sheaf $\LL$ on $X$, such that for
all $\alpha > 0$ the sheaf $S^{\alpha \beta}(\GG) \otimes \LL^{\beta}$
is globally generated over $U$ for some $\beta > 0$.
\end{theorem}

\proof $(a) \Leftrightarrow (b)$ is just the definition of weak 
positivity, and we obviously have
$(e) \Rightarrow (c) \Rightarrow (b) \Rightarrow  (d) \Rightarrow
(f) $. So it remains to show $(f) \Rightarrow (e)$. Let $\HH$ be
ample with respect to $U$ and let $\alpha > 0$. By 
\ref{lokampkrit} we find 
a blowing up $\tau : X' \to X$, an ample sheaf $\HH'$ on
 $X'$ and some $\mu > 0$ with an inclusion
 $\HH' \to \tau^{\ast}\HH^{\mu}$, being an isomorphism over
$\tau^{-1}(U)$.
 Let $\LL' = \tau^{\ast}\LL$ and let us choose $\gamma > 0$  such that
 ${\LL'}^{-1} \otimes {\HH'}^{\gamma}$ is globally generated. Thus,
 for some $r > 0$, we 
have a surjective map $\bigoplus^r \LL' \to {\HH'}^{\gamma}$.
By assumption, there exists some $\beta' > 0$ such that, for some 
$s > 0$, we find a map
 $$\bigoplus^s \OO_{X'} \to \tau^{\ast} \left( S^{2 \alpha \gamma \mu
\beta'}
 (\GG) \otimes \LL^{\beta'} \right)
 =  \tau^{\ast} S^{2 \alpha \gamma \mu \beta'}
 (\GG) \otimes {\LL'}^{\beta'},$$
 surjective over $\tau^{-1}(U)$.
 Thus we obtain maps
$$\bigoplus^r \bigoplus^s \OO_{X'} \to \bigoplus^r \left( \tau^{\ast}
 S^{2 \alpha \gamma \mu \beta'} (\GG) \otimes {\LL'}^{\beta'} \right)
 = \tau^{\ast}  S^{2 \alpha \gamma \mu \beta'} (\GG) \otimes
 \bigoplus^r {\LL'}^{\beta'}$$
and hence $$\bigoplus^r \bigoplus^s \OO_{X'} \to \tau^{\ast}  S^{2
\alpha \gamma \mu \beta'} (\GG) \otimes {\HH'}^{\gamma \beta'} \to
\tau^{\ast} \left( S^{2 \alpha \gamma \mu \beta'} (\GG) \otimes
\HH^{\mu \gamma \beta'} \right),$$ all being surjective over
$\tau^{-1}(U)$. This induces a map $$\bigoplus \tau_{\ast}\OO_{X'}
\otimes \HH^{\mu \gamma \beta'} \to S^{2 \alpha \gamma \mu \beta'}
(\GG) \otimes \HH^{2 \mu \gamma \beta'},$$ surjective over $U$.
We may assume that the sheaf on the left hand side is globally
generated, and, taking $\beta = 2 \gamma \mu \beta'$, this implies $(e)$.
\bewend

We consider the following class of ideal sheaves: 

\begin{definition}\label{ivoll}
\rm Let $X$ be a normal quasi-projective variety and let $\II$ be an 
ideal sheaf on $X$.
Let $\tau : X' \to X$ be a birational morphism such that $X'$ is
normal and $\II' = \tau^{-1}\II \cdot \OO_{X'}$ is invertible. 

1. We call $\II$ {\it full}, if the natural map $\II
\to \tau_{\ast}\II'$ is an isomorphism.

2. If $\LL$ is an invertible sheaf on $X$ and 
$U \subseteq X$ an open subset, then we call $\LL \otimes
\II$ {\it weakly positive over $U$}, if $\tau^{\ast}\LL \otimes
\II'$ is weakly positive over $\tau^{-1}(U)$.

3. If $X$ is projective, then we call $\LL \otimes \II$ {\it numerically
effective (nef)}, if $\tau^{\ast}\LL \otimes \II'$ is numerically
effective on $X'$.
\end{definition}

\begin{lemma}\label{ivoll3}
Let $X$ be a normal quasi-projective  variety, let $\LL$ be an invertible
sheaf on $X$, let $\II$ be a full ideal sheaf on $X$ and let 
$U \subseteq X$ be an open subset. 

(a) $\LL \otimes \II$ is weakly positive over
$U$ if and only if for every sheaf $\HH$, ample with respect to $U$, and 
for all $\alpha
> 0$ the sheaf $\LL^{\alpha \beta} \otimes \II^{\alpha \beta}
\otimes \HH^{\beta}$ is globally generated over $U$ for some $\beta > 0$.

(b) If $\LL_1$ and $\LL_2$ are invertible sheaves on $X$ and
if for every $\mu > 0$ there is some $\nu> 0$ such that $\LL^{\mu
\nu} \otimes \II^{\mu \nu} \otimes \LL_1^{\nu} \otimes \LL_2$ is
weakly positive over $U$, then so is $\LL \otimes \II$.
\end{lemma}

\proof $(b)$ follows immediately from $(a)$ and from \cite{ev84},
4.3. In order to show $(a)$, let us assume $\LL \otimes \II$ is weakly 
positive over
$U$. We choose some birational morphism $\tau : X' \to X$ such
that $X'$ is normal, $\II' = \tau^{-1} \II \cdot \OO_{X'}$ is
invertible on $X'$ and $\tau|_{\tau^{-1}(U)}$ is an isomorphism.
If $\HH$ is ample with respect to $U$, then $\tau^{\ast}\HH$ is ample 
with respect to $\tau^{-1}(U)$ by
\ref{lokampend}.
If we choose $\alpha > 0$, then by \ref{amplunabh} there exists some
$\beta > 0$ such that
$(\tau^{\ast}\LL \otimes \II')^{2 \alpha \beta} \otimes
\tau^{\ast}\HH^{\beta}$ is globally generated over $\tau^{-1}(U)$.
Thus there exists a map

$$\bigoplus \OO_{X'} \longrightarrow (\tau^{\ast}\LL \otimes \II')^{2
\alpha \beta} \otimes \tau^{\ast}\HH^{\beta},$$ surjective over
$\tau^{-1}(U)$, and hence we obtain a map
\begin{eqnarray*}
\bigoplus \tau_{\ast}\OO_{X'} \otimes \HH^{\beta} & \longrightarrow &
\tau_{\ast} \left((\tau^{\ast}\LL \otimes \II')^{2 \alpha \beta} \otimes
\tau^{\ast}
\HH^{\beta}\right) \otimes \HH^{\beta} \\
& & = \LL^{2 \alpha \beta} \otimes \tau_{\ast}{\II'}^{2 \alpha \beta}
\otimes \HH^{2 \beta} = \LL^{2 \alpha \beta} \otimes \II^{2 \alpha \beta}
\otimes
\HH^{2 \beta},
\end{eqnarray*}
surjective over $U$, where the last equality holds because $\II$
is full. Now we may assume that the sheaf on the left hand side is 
globally generated, which proves the necessity of the condition.

In order to show it is sufficient, too, let us choose again some
birational
morphismus $\tau : X' \to X$ with the properties from above. Let
$\alpha > 0$ and let $\HH$ be ample with respect to $U$, hence 
$\tau^{\ast}\HH$ is ample with respect to $\tau^{-1}(U)$ by 
\ref{lokampend}. Then $\tau^{\ast}\left(\LL^{\alpha
\beta} \otimes \II^{\alpha \beta} \otimes \HH^{\beta}\right)$ is
globally generated over $\tau^{-1}(U)$. By definition there exists
a map $$\tau^{\ast}\left(\LL^{\alpha \beta} \otimes \II^{\alpha
\beta} \otimes \HH^{\beta}\right) \to (\tau^{\ast}\LL \otimes
\II')^{\alpha \beta} \otimes \tau^{\ast}\HH^{\beta}, $$ surjective
over over $\tau^{-1}(U)$. This implies the weak positivity of
$\tau^{\ast}\LL \otimes \II'$ over $\tau^{-1}(U)$ and so by
definition the weak posititity of $\LL \otimes \II$ over $U$, which
completes the proof of $(a)$.
\bewend

We shall need two more properties:

\begin{lemma}\label{posrunter}
Let $\varrho : Y \to X$ be a surjective morphism of normal
quasi-projective  varieties and let $U \subseteq X$ be an open subset such
that $\varrho|_{\varrho^{-1} (U)}$ is finite. Let $\GG$ be a
locally free sheaf on $Y$, let $\LL$ be an invertible sheaf and 
let $\JJ$ be an
ideal sheaf on $X$ such that for every $l \in \N-\{0\}$ the sheaf
$\JJ^l$ is full and there exists a map $$\varrho_{\ast}S^l(\GG)
\longrightarrow \LL^l \otimes \JJ^l,$$ surjective over $U$. If in
this situation $\GG$ is weakly positive over $\varrho^{-1}(U)$,
then $\LL \otimes \JJ$ is weakly positive over $U$.
\end{lemma}

\proof
Let $\alpha > 0$ and let $\HH$ be an ample invertible sheaf on $X$.
Since $\varrho$ is finite over $U$, $\varrho^{\ast}\HH$ is ample
over $\varrho^{-1}(U)$ by \ref{lokampend}. If $\GG$ is weakly
positive over $\varrho^{-1}(U)$, then by \ref{amplunabh} there
exists some $\beta
> 0$ such that the sheaf $S^{ 2 \alpha \beta}(\GG) \otimes
\varrho^{\ast}\HH^{\beta}$ is globally generated over
$\varrho^{-1}(U)$. Thus we obtain a map $$\bigoplus \OO_Y \to S^{ 2
\alpha \beta}(\GG) \otimes \varrho^{\ast} \HH^{\beta}, $$
surjective over $\varrho^{-1}(U)$ and hence a map $$\bigoplus
\varrho_{\ast}\OO_Y \otimes  \HH^{\beta} \to \varrho_{\ast} \left(
S^{ 2 \alpha \beta}(\GG) \otimes \varrho^{\ast}\HH^{\beta} \right)
\otimes \HH^{\beta} = \varrho_{\ast}S^{ 2 \alpha \beta}(\GG)
\otimes \HH^{2 \beta},$$ surjective over $U$. By assumption and
since $\JJ^l$ is full for every $l \in \N - \{0\}$, there exists a
map $$\bigoplus \varrho_{\ast}\OO_Y \otimes  \HH^{\beta} \to
\LL^{2 \alpha \beta} \otimes \JJ^{2 \alpha \beta} \otimes \HH^{2
\beta},$$ surjective over $U$.
We may assume that the sheaf on the left hand side is globally generated,
therefore $\LL^{2 \alpha \beta} \otimes \JJ^{2 \alpha \beta}
\otimes \HH^{2 \beta}$
is globally generated over $U$ which by \ref{ivoll3} yields the weak
positivity of $\LL \otimes \II$ over $U$. \bewend

\begin{lemma}\label{ueberh}
Let $X$ be a normal quasi-projective variety, let $\GG$ be a locally free
sheaf on $X$ and let $U\subseteq X$ be an open subset. Then the
following statements are equivalent:

(a) $\GG$ is weakly positive over $U$.

(b) There exists some $\mu > 0$ and some sheaf $\HH$ on $X$, ample 
with respect to $U$, such
that for every morphism $\tau : X' \to X$, finite over $U$, 
with $\tau^{\ast}\HH = {\HH'}^{\delta}$ for some $\delta
> 0$ and some invertible sheaf $\HH'$ on $X'$, the sheaf
$\tau^{\ast}\GG \otimes {\HH'}^{\mu}$ is weakly positive over
$\tau^{-1}(U)$.

(c) There exists some $\mu > 0$, some invertible sheaf $\LL$ on
$X$ and for every $\delta \in \N-\{0\}$ a morphism $\tau_{\delta}
: X^{\delta} \to X$, finite over $U$, such that
$\tau_{\delta}^{\ast}\LL = {\LL'}^{\delta}$ for some invertible sheaf
$\LL'$ on $X^{\delta}$ and such that $\tau_{\delta}^{\ast}\GG
\otimes {\LL'}^{\mu}$ is weakly positive over
$\tau_{\delta}^{-1}(U)$.

\end{lemma}

\proof We only have to show $(c) \Rightarrow (a)$. Let $\alpha > 0$
and let $\HH$ be an 
ample invertible sheaf on $X$ such that $\LL
\otimes \HH$ is ample, too. By
assumption there exists a morphism $\tau : X' \to X$ , finite over
$U$, such that $\tau^{\ast}\LL = {\LL'}^{1 + 2 \alpha
\mu}$ for some invertible sheaf $\LL'$ on $X'$. By \cite{v95}, 2.1 we
may assume that in addition $\tau^{\ast}\HH = {\HH'}^{1 + 2
\alpha \mu}$ for some invertible sheaf $\HH'$ on $X'$. Since 
$\tau^{\ast}\GG \otimes {\LL'}^{\mu}$ is weakly positive over
$\tau^{-1}(U)$, so is $\tau^{\ast}\GG
\otimes (\LL' \otimes \HH')^{\mu}$. Hence, for some
$\beta > 0$ we obtain a map
$$\bigoplus \OO_{X'} \to S^{2 \alpha \beta'}\left(\tau^{\ast}\GG
\otimes (\LL' \otimes \HH')^{\mu}\right) \otimes (\LL' \otimes
\HH')^{\beta'} = $$ $$\tau^{\ast}\left(S^{2 \alpha
\beta'}(\GG)\right) \otimes (\LL' \otimes \HH')^{(1 + 2 \alpha \mu)
\beta'}, $$ surjective over $\tau^{-1}(U)$ and thus a map
$$\bigoplus \tau_{\ast}\OO_{X'} \otimes (\LL \otimes \HH)^{\beta'}
\to S^{2 \alpha \beta'}(\GG) \otimes (\LL \otimes \HH)^{2
\beta'},$$ surjective over $U$. We may assume that the sheaf on the
left hand side
is globally generated and we find that $S^{2 \alpha \beta'}(\GG)
\otimes (\LL \otimes \HH)^{2 \beta'}$ is globally generated over
$U$. Choosing $\beta = 2 \beta'$, we are done. \bewend

Finally, let us recall the Fujita-Kawamata Positivity Theorem in a
slightly modified version. A proof due to Koll\'ar can be found in
\cite{v95}, 2.41.

\begin{theorem}[Fujita \cite{f78}, Kawamata \cite{k81}]\label{posfk}
Let $f : X \to Y$ be a projective surjective morphism of smooth
quasi-projective varieties and let $Y_0 \subseteq Y$ be an open subset
such that $f|_{f^{-1}(Y_0)}$ is smooth and $f_{\ast}\omega_{X/Y}$ is
locally free
over $Y_0$. Then $f_{\ast}\omega_{X/Y}$ is weakly positive over $Y_0$.
\end{theorem}

\begin{corollary}\label{posfkkor}
In the situation of \ref{posfk}, let $\LL$ be an invertible sheaf
on $X$ and let $D$ be an effective normal crossing divisor on $X$ such
that $\LL^N(-D)$ is semi-ample. Then $f_{\ast}(\LL \otimes
\noidealtext{X/Y}{D}{N})$ is weakly positive over $Y_0$.
\end{corollary}

\proof For $\beta > 0$ sufficiently large, we have $\LL^{\beta
\cdot N} = \OO_X(\Gamma + \beta \cdot D)$, where $\Gamma$ is a
smooth divisor and $\Gamma + \beta \cdot D$ is a normal crossing
divisor, and on the other hand we have $\noidealtext{X/Y}{\Gamma +
\beta \cdot D}{\beta \cdot N} = \noidealtext{X/Y}{D}{N}$. We thus
may assume that $\LL^N = \OO_X(D)$, and the statement follows
from considering the cyclic covering according to this situation
(see \cite{v95}, 2.43). \bewend

\begin{center}
\section{\sc Some Positivity Statements}
\end{center}

\begin{definition}\label{defmult}
\rm Let $X$ be a normal variety with at most rational singularities
and let $\Gamma$ be an effective Cartier divisor
on $X$. Let $\tau : X' \to X$ be a blowing up such that both $X'$
is smooth and $\Gamma' = \tau^{\ast}\Gamma$ has normal crossings.
For every $N \in \N-\{0\}$ let us define $$\oideal{X}{\Gamma}{N} =
\tau_{\ast}\omega_{X'} \left( - \left[\frac{\Gamma'}{N} \right]
\right).$$ 
Due to \cite{v95}, 5.10, this definition does not depend on the chosen
blowing up.
Moreover, we define 
$$\CC_X(\Gamma,N) = \coker \left(
\oideal{X}{\Gamma}{N} \longrightarrow \omega_X \right).$$ 
Since $\omega_X =
\tau_{\ast}\omega_{X'}$, we have $\CC_X(\Gamma,N) = 0$ for $N$
sufficiently large, so we define $$e(\Gamma) = \min \Big\{ N
\in \N-\{0\} \, ; \, \CC_X(\Gamma,N) = 0 \Big\}.$$ If in addition 
$X$ is projective and $\LL$ is an effective invertible sheaf on
$X$, then let $$e(\LL) = \max \Big\{ e(\Gamma) \, ; \, \Gamma
\mbox{\,  effective Cartier divisor with \,} \LL = \OO_X(\Gamma)
\Big\}.$$

\end{definition}

\begin{remark}\label{multbir2}
Immediately from the definition we see that, given any 
birational morphism $\delta : Z \to X$ of varieties and any 
effective Cartier divisor $\Gamma$ on $X$, we may choose
a blowing up $\tau : X' \to
X$ in such a way that $X'$ is smooth and $\Gamma' =
\tau^{\ast}\Gamma$ is a normal crossing divisor on $X'$ and that
in addition $\tau$ factors over $\delta$, hence
$$\delta_{\ast}\oideal{Z}{\delta^{\ast}\Gamma}{N} =
\oideal{X}{\Gamma}{N}.$$
\end{remark}

\begin {lemma}\label{multetal}
If $\sigma : Y \to X$ is a finite surjective morphism of smooth
varieties,
 $\Gamma$ an effective Cartier divisor
 on $X$ and $N \in \N-\{0\}$, then there exists a map
 $$\oideal{Y}{\sigma^{\ast}\Gamma}{N} \longrightarrow \sigma^{\ast}
\oideal{X}{\Gamma}{N}.$$
If $U \subseteq X$ is the open subset such that
$\sigma|_{\sigma^{-1}(U)}$ is \'etale, then this map is an
isomorphism over $\sigma^{-1}(U)$.
 \end{lemma}

 \proof
Let $\tau : X'\to X$ be a blowing up such that $X'$ is smooth and
$\Gamma' =
 \tau^{\ast}\Gamma$ is a normal crossing divisor.
Let $Y'$ be a desingularization of the fibre product $Y \times_X
X'$ such that we have the following diagram:
  \begin{eqnarray*}
 Y' & \stackrel{\sigma'}{\longrightarrow} & X' \\
 \tau' \downarrow & & \downarrow \tau \\
 Y & \stackrel{\sigma}{\longrightarrow} & X .
 \end{eqnarray*}
We may assume that ${\sigma'}^{\ast}\Gamma' = {\tau'}^{\ast}
\sigma^{\ast}\Gamma$ is a normal crossing divisor on $Y'$. We have
$$\oideal{Y}{\sigma^{\ast}\Gamma}{N} = \tau'_{\ast}\noideal{Y'}
{{\tau'}^{\ast}\sigma^{\ast}\Gamma}{N} = \tau'_{\ast}\noideal{Y'}
{{\sigma'}^{\ast}\Gamma'}{N}.$$
Considering integral parts of $\Q$-divisors, there exists an 
injective map
$$\noideal{Y'}{{\sigma'}^{\ast}\Gamma'}{N} \to {\sigma'}^{\ast}
\noideal{X'}{\Gamma'}{N},$$ which is an isomorphism over
$\sigma^{-1}(U)$, and in this case even
the multiplicities of the components do not change. We obtain a map
$$\oideal{Y}{\sigma^{\ast}\Gamma}{N} \longrightarrow \tau'_{\ast}
{\sigma'}^{\ast} \noideal{X'}{\Gamma'}{N}.$$ Since $\sigma$ is flat
by \cite{h77}, III. Exercise 9.3, we obtain by flat base change
(\cite{h77}, III.9.3) a map 
 $$\oideal{Y}{\sigma^{\ast}\Gamma}{N} \longrightarrow 
\tau'_{\ast}{\sigma'}^{\ast}
\noideal{X'}{\Gamma'}{N} = \sigma^{\ast}\tau_{\ast}
 \noideal{X'}{\Gamma'}{N} = \sigma^{\ast} \oideal{X}{\Gamma}{N}, $$
being an isomorphism over $\sigma^{-1}(U)$. \bewend

We shall need the following fact which generalizes \cite{v95}, 5.21:
\begin{proposition}\label{products}
Let $X_1 , \dots , X_r$ be smooth projective varieties and let $X =
X_1 \times \dots \times X_r$ denote their product. For $i = 1 , \dots
, r$ let us consider effective divisors $D_i$ and effective 
invertible sheaves $\LL_i$ on $X_i$. Let $\Gamma$ denote the 
induced divisor on $X$ and let $\LL$ denote the induced 
invertible sheaf on $X$. Then we have

(a) $e(\Gamma) = \max \{e(D_1) , \dots , e(D_r)\}$.

(b) $e(\LL) = \max \{e(\LL_1) , \dots , e(\LL_r)\}$.
\end{proposition}

\proof 
We may restrict ourselves to the case of two factors. In order to show
$(a)$, we may assume
$e = e (D_2) \geq
e(D_1)$. 
Let us start with the case where $D_2$ is a normal crossing divisor on
$X_2$. Let $p_2 : X_1 \times X_2 \to X_2$ be the second
projection. For $\Gamma = D_1 \times X_2$ one has
$e(\Gamma|_{p_2^{-1}(x)}) = e (D_1)$ for all fibres of the second
projection, and by \cite{v95}, 5.18 one has
 $$\oideal{X_1 \times X_2}{D_1 \boxplus D_2}{N}  =
\omega_{X_1 \times X_2}\left(-p_2^{\ast}\left[\frac{D_2}{N}\right]
\right)$$ for $N \geq e(D_1)$. In particular this holds for $N=e$, 
and since in this case the sheaf on the right hand side is equal to 
$\omega_{X_1 \times X_2}$, we have 
$e(D_1 \boxplus D_2) \leq e$.

Let now $D_2$ be an arbitrary effective divisor on $X_2$ and $N
\geq e(D_1)$. We consider a blowing up $\tau : X_2' \to X_2$ such
that $X_2'$ is smooth and $D_2' = \tau^{\ast}D_2$ is a normal
crossing divisor on $X_2'$:
\begin{eqnarray*}
X_1 \times X_2' & \stackrel{\tau'}{\longrightarrow} & X_1 \times X_2 \\
\downarrow p_2' & & \downarrow p_2 \\
X_2' & \stackrel{\tau}{\longrightarrow} &  X_2.
\end{eqnarray*}
For $D_1 \boxplus D_2'$ on $X_1 \times X_2'$ we use the
first case and the compatibility of relatively canonical sheaf with 
base change to obtain
\begin{eqnarray*}
\oideal{X_1 \times X_2'}{D_1 \boxplus D_2'}{N} & = &
\omega_{X_1 \times X_2'}\left(-{p_2'}^{\ast}\left[\frac{D_2'}{N}\right]
\right) \\
& = & \omega_{(X_1 \times X_2')/X_2'} \otimes
{p_2'}^{\ast}\noideal{X_2'}{D_2'}{N} \\
& = & {\tau'}^{\ast}\omega_{(X_1 \times X_2)/X_2} \otimes
{p_2'}^{\ast}\noideal{X_2'}{D_2'}{N}.
\end{eqnarray*}  

Using \ref{multbir2} and flat base change (\cite{h77}, III.9.3), 
this yields 
$$\oideal{X_1 \times X_2}{D_1 \boxplus
D_2}{N} = \omega_{(X_1 \times X_2)/X_2} \otimes p_2^{\ast}
\tau_{\ast}\noideal{X_2'}{D_2'}{N}.$$ 
For $N = e$ the right hand
side is nothing but $\omega_{(X_1 \times X_2)/X_2} \otimes
p_2^{\ast} \omega_{X_2} = \omega_{X_1 \times X_2}$, and we obtain
$e(D_1 \boxplus D_2) \leq e$ for this case, too.

To show that $e$ is a lower bound, too, it suffices to show that 
for an open subset $U \subseteq X_1 \times X_2$ one has 
$e\left((D_1 \boxplus
D_2)|_U\right) \geq e$. We thus may assume that $D_1 = 0$, and in
the above calculation may choose $N = e-1$. Then, since $e(D_2) =
e$, we obtain that $$\oideal{X_1 \times X_2}{D_1 \boxplus
D_2}{e-1} \longrightarrow \omega_{X_1 \times X_2}$$ is not an 
isomorphism, so $e-1 < e(D_1 \boxplus D_2)$. 

For $(b)$ we may assume $e(\LL_2) \geq e(\LL_1)$. Let us
choose effective divisors $D_1$ on $X_1$ and $D_2$ on $X_2$ with
$e(\LL_1) = e(D_1)$ and $e(\LL_2) = e(D_2)$. Then $D_1 \boxplus
D_2$ is a section of $\LL$, and by $(a)$ we obtain $$e(\LL)
\geq e(D_1 \boxplus D_2) = e(D_2) = e(\LL_2).$$

For the other direction we consider the projections $p_1 : X \to
X_1$ and $p_2 : X \to X_2$. One has $e(\LL|_{p_1^{-1}(x_1)}) =
e(\LL_2)$ for every fibre $p_1^{-1}(x_1)$ and
$e(\LL|_{p_2^{-1}(x_2)}) = e(\LL_1)$ for every fibre
$p_2^{-1}(x_2)$. Let $\Gamma$ be an effective divisor with $\LL =
\OO_X(\Gamma)$. By \cite{v95}, 5.19, the support of
$\CC_X\left(\Gamma,e(\LL_2)\right)$ is of the form $p_1^{-1}(S_1) =
S_1 \times X_2$ for some closed subvariety $S_1 \subseteq X_1$ and
the support of $\CC_X\left(\Gamma,e(\LL_1)\right)$ is of the form
$p_2^{-1}(S_2) = X_1 \times S_2$ for some closed subvariety $S_2
\subseteq X_2$. But the support of $\CC_X\left(\Gamma,e(\LL_2)\right)$ 
is contained in the support
of $\CC_X\left(\Gamma,e(\LL_1)\right)$,
which yields the vanishing of $\CC_X\left(\Gamma,e(\LL_2)\right)$ or 
in other words $e(\LL)
\leq  e(\LL_2)$. \bewend

The usual vanishing and positivity theorems (see \cite{v95}, for
example) can be extended to
the situation described above:

\begin{theorem}\label{versmult}
Let $X$ be a normal projective variety with at most rational
singularities, let $\LL$ be an invertible sheaf on $X$, and let 
$D$ be an effective Cartier divisor on $X$ and $N \in \N-\{0\}$.

(a) If $\LL^N(-D)$ is nef and big, then for all $i > 0$ one has
$$H^i\left(X,\LL \otimes \oideal{X}{D}{N} \right) = 0.$$

(b) If $\LL^N(-D)$ is semi-ample and $(\LL^N(-D))^{\nu}(-B)$ is effective
for some
effective Car\-tier divisor $B$ on $X$ and some $\nu \in
\N-\{0\}$, then for all $i > 0$ the map $$H^i \left(X, \LL(B)
\otimes \oideal{X}{D}{N} \right) \longrightarrow H^i \left(B,\left(
\LL(B) \otimes \oideal{X}{D}{N} \right)|_B \right)$$

is surjective.

(c) Let $f : X \to Y$ be a projective surjective morphism on a
smooth quasi-projective variety $Y$ and $Y_0 \subseteq Y$ the open
subset such that $f_{\ast}(\LL \otimes \oidealtext{X/Y}{D}{N})$ is
locally free over $Y_0$ and $f|_{f^{-1}(Y_0)}$ is smooth. If
$\LL^N(-D)$ is semi-ample, then $f_{\ast}(\LL \otimes
\oidealtext{X/Y}{D}{N})$ is weakly positive over $Y_0$.
\end{theorem}

\proof Let $\tau : X' \to X$ a blowing up such that $D' =
\tau^{\ast}D$ is a normal crossing divisor and $X'$ is smooth.
Writing $\LL' = \tau^{\ast}\LL$, the sheaf ${\LL'}^N(-D')$ is nef
and big on $X'$, and writing $B' = \tau^{\ast}B$, the sheaf
$({\LL'}^N(-D'))^{\nu}(-B')$ is effective. 
We obtain $(a)$ and $(b)$ from 
the corresponding vanishing theorems for integral parts of
$\Q$-divisors which hold on $X'$ (see \cite{v95}, 2.28 and 2.33, 
respectively). Finally, $(c)$
follows from \ref{posfkkor}: If $f' : X' \to Y$ is the induced
map, then $$f'_{\ast}\left( \LL' \otimes
\noideal{X'/Y}{D'}{N}\right) = f_{\ast} \left( \LL \otimes
\oideal{X/Y}{D}{N}\right)$$ is weakly positive over $Y_0$. \bewend

The following two statements correspond to \cite{ev84}, 6.5 and 6.4,
respectively. 

\begin{lemma}\label{surj}
Let $Z$ and $B$ be smooth quasi-projective varieties and let 
$p : Z \to B$ be a projective surjective morphism of relative 
dimension $k$. Let us consider the open subset $W' \subseteq B$ 
such that $p|_{p^{-1}(W')} : p^{-1}(W') \to W'$ is smooth.
Let $\MM$ be an invertible sheaf on $Z$,
let $D$ be an effective divisor on $Z$ and let $U \subseteq Z$ be an
open subset such that  
$\MM^N(-D)$ is relatively semi-ample over $U$ and 
$e(D|_{Z_b}) \leq N$ for all $b \in W'$.
Let $\KK$ be an invertible sheaf on $Z$ such that 
$\KK|_{Z_b}$ is ample with respect to $Z_b \cap U$
and such that $\varepsilon
= \varepsilon (\KK|_{Z_b},Z_b \cap U) \geq 1$ for all $b \in W'$ (here
$\varepsilon$ denotes the Seshadri index, see \cite{el93}, for example). 
Then there exists a nonempty open subset $W \subseteq B$ such that
$$p^{\ast}p_{\ast}\left( \MM \otimes \KK^{k+1} \otimes
\oideal{Z/B}{D}{N} \right) \to \MM \otimes \KK^{k+1} \otimes
\oideal{Z/B}{D}{N}$$ over $U \cap p^{-1}(W)$ is surjective.
\end{lemma}

\proof
By \cite{v95}, 5.10. we may choose 
an open subset $W \subseteq W'$ where 
$\oideal{Z/B}{D}{N}|_{Z_b} =
\oideal{Z_b}{D|_{Z_b}}{N}$ for all $b \in W$. (In fact, considering 
some blowing up $\tau' : Z'
\to Z$ such that $D' = \tau^{\ast}D$ is a normal crossing divisor, we
may take $W$ as the set of all $b \in W$ such that $D'|_{Z'_{b}}$ is a
normal crossing divisor.)

Let $b \in W$ be in general position and let us fix some $z 
\in Z_b \cap U$. 
From now on, the index $b$ means restriction to the fibre $Z_b$.
We 
consider the blowing up $\tau : \tilde{Z_b} \to Z_b$ of the
fibre $Z_b$ in $z$ and we denote by $E_z$ the exceptional divisor
of $\tau$. Then, by \cite{h77}, II. Exercise 8.5, we obtain that
$\omega_{\tilde{Z_b}} = \tau^{\ast}\omega_{Z_b} \otimes
\OO_{\tilde{Z_b}} \left((k-1) \cdot E_z\right)$. Moreover, let $\varrho
: Z_b' \to Z_b$ be a blowing up such that ${\varrho}^{\ast}D_b$ 
is a
normal crossing divisor and $Z_b'$ is smooth. We may assume that
$\varrho$ factors as $\varrho = \tau \circ \eta$ for some morphism
$\eta : Z_b' \to\tilde{Z_b}$. We write 
$\omega_{Z_b'} = {\varrho}^{\ast}\omega_{Z_b} \otimes
\OO_{Z_b'}(E+F)$, where $E$ is the part of the exceptional divisor
with $\varrho(E) = z$.
We now have to prove that the sheaf
$\MM_b \otimes \KK_b^{k+1} \otimes
\oidealtext{Z_b} {D_b}{N}$ is globally generated over $Z_b \cap U$. 

We are obviously allowed to replace $N$
and $D$ by multiples and, hence, may assume that $\MM^N(-D)$ is
relatively globally generated over $U$. Considering the subsheaf
of $\MM^N(-D)$ which is globally generated and replacing $Z$ by a
blowing up with centre outside of $U$, making this sheaf
invertible, we may assume that $\MM^N(-D)$ is
relatively globally generated, hence relatively numerically effective.

By assumption, $\tau^{\ast}\MM_b^N(-D_b)$ and
$\tau^{\ast}\KK_b^{N \cdot k} \otimes \OO_{\tilde{Z_b}}(-N \cdot k
\cdot E_z)$ are nef, and $\tau^{\ast}\KK_b^N$ is nef and big, hence 
$\tau^{\ast}(\MM_b \otimes
\KK_b^{k+1})^N \otimes \OO_{\tilde{Z_b}}(-N \cdot k \cdot E_z -
\tau^{\ast}D_b)$ is nef and big. By \cite{v95}, 5.22 we obtain
a surjection
$$H^0 \left( Z_b , \MM_b \otimes \KK_b^{k+1} \otimes
\omega_{Z_b} \right) \to H^0 \left( Z_b , 
\CC_b\right),$$
where $\CC_b$ denotes the cokernel of 
$$\MM_b \otimes \KK_b^{k+1} \otimes
{\varrho}_{\ast}\omega_{Z_b'}\left( - \left[ \frac{{\varrho}^{\ast}
D_b}{N} \right] -
k \cdot \eta^{\ast}E_z \right) \to \MM_b \otimes \KK_b^{k+1}.$$
By assumption, the support of ${\CC_b}_{|U \cap Z_b}$ is contained 
in \{z\}.
Assuming ${\CC_b}_{|U \cap Z_b} = 0$, we obtain 
$${\varrho}_{\ast}\omega_{Z_b'}\left( - \left[
\frac{{\varrho}^{\ast}D_b}{N} \right] - k \cdot \eta^{\ast}E_z \right) =
{\varrho}_{\ast}\omega_{Z_b'}\left( - \left[ \frac{{\varrho}^{\ast}D_b}{N}
\right] \right) = \oideal{Z_b}{D_b}{N} = \omega_{Z_b},$$ and 
hence
\begin{center}
\begin{tabular}{l}
${\varrho}_{\ast} \left( {\varrho}^{\ast}\omega_{Z_b} \otimes \OO_{Z_b'}
\left(E - \left[ \frac{{\varrho}^{\ast}D_b}{N}
\right] - k \cdot \eta^{\ast}E_z \right) \right)$ \\ \\
$= \omega_{Z_b} \otimes {\varrho}_{\ast}
\OO_{Z_b'}
\left(E - \left[ \frac{{\varrho}^{\ast}D_b}{N}
\right] - k \cdot \eta^{\ast}E_z \right) = \omega_{Z_b}.$
\end{tabular}
\end{center}
But this means $$E \geq \left[ \frac{{\varrho}^{\ast}D_b}{N} \right] + k
\cdot \eta^{\ast}E_z, $$ so $E$ 
contains $\eta^{\ast}E_z$ with a multiplicity of at least $k$, which
contradicts the definition of $E$. So $\CC_b$ is concentrated in the
point $\{z\}$, and we may apply \cite{k99}, 2.11 to obtain that 
$\MM_b \otimes \KK_b^{k+1} \otimes
\oidealtext{Z_b} {D_b}{N}$ is globally generated over $Z_b \cap
U$. \bewend

\begin{proposition}\label{prop}
Keeping the assumptions from \ref{surj}, let 
$\LL$ and $\BB$ be invertible sheaves on $Z$ such that $\LL(-D)$ 
is semi-ample over $U$, $\LL$ 
is relatively semi-ample over $U$ and $\BB$ is
relatively numerically effective. Let us assume moreover that 
$e(D|_U) \leq N$ and 
that, for some $r > 0$, there exists a map 
$\bigoplus^r \omega_{Z/B} \to \BB$, surjective over $U$. 
Then there exists a nonempty open subset $W \subseteq B$ such that
$\LL \otimes \BB^N \otimes \KK^{k+1}$ is weakly positive over $U
\cap p^{-1}(W)$.
\end{proposition}

\proof
Let us start with the case where $p : Z \to B$ is flat. 
Taking $\MM = \LL \otimes \BB^{N-1}$, we see that 
$\MM^N(-D)$ is relatively semi-ample over $U$. Using \ref{surj} 
and since $e(D|_U) \leq N$, we obtain maps
\small
\begin{eqnarray*}
p^{\ast}p_{\ast}\left( \LL \otimes \BB^{N-1} \otimes \KK^{k+1} \otimes
\oideal{Z/B}{D}{N} \right) &  \to &  \LL \otimes \BB^{N-1} \otimes
\KK^{k+1}
\otimes \oideal{Z/B}{D}{N} \\
&  \to & \LL \otimes \BB^{N-1} \otimes \KK^{k+1}
\otimes \omega_{Z/B},
\end{eqnarray*}
\normalsize 
surjective over $U \cap p^{-1}(W'')$
for some nonempty open subset $W'' \subseteq B$. If
$\tilde{\HH}$ is an ample invertible sheaf on $B$, then the sheaf
$$p_{\ast}\left(\LL \otimes \BB^{N-1} \otimes \KK^{k+1} \otimes
\oideal{Z/B}{D}{N}\right) \otimes \tilde{\HH}^{\mu}$$ is globally
generated for $\mu$ sufficiently large, and with $\HH =
p^{\ast}\tilde {\HH}$ we obtain that $$p^{\ast}p_{\ast}\left(\LL
\otimes \BB^{N-1} \otimes \KK^{k+1} \otimes
\oideal{Z/B}{D}{N}\right) \otimes \HH^{\mu}$$ is globally
generated, too. Since by \cite{v95}, 2.16 quotient sheaves inherit weak
positivity, we obtain that $$\LL \otimes
\BB^{N-1} \otimes \KK^{k+1} \otimes \omega_{Z/B} \otimes
\HH^{\mu}$$ is weakly positive over $U \cap p^{-1}(W'')$. By
assumption we have a map
\begin{eqnarray*} \bigoplus^r
\left(\LL \otimes \BB^{N-1} \otimes \KK^{k+1} \otimes \omega_{Z/B}
\otimes \HH^{\mu} \right) &  = \LL \otimes \BB^{N-1} \otimes
\KK^{k+1} \otimes \bigoplus^r \omega_{Z/B} \otimes \HH^{\mu} & \\
& \to \LL \otimes \BB^N \otimes \KK^{k+1} \otimes \HH^{\mu},&
\end{eqnarray*}
surjective over $U$, which yields, using \cite{v95}, 2.16 again,
the weak positivity of $\LL \otimes \BB^N \otimes \KK^{k+1}
\otimes \HH^{\mu}$ over $U \cap p^{-1}(W'')$, where $\mu$
is sufficiently large. This gives sense to the following
definition: $$
\begin{array}{cccc}
\gamma & = & \min & \Big\{   \mu \in N \cdot \N \, ; \, \LL
\otimes \BB^N \otimes \KK^{k+1} \otimes \HH^{\mu} \mbox {\, weakly
positive over \,} \\ & & & U \cap p^{-1}(W) \mbox {\, for a
nonempty open subset \,} W \subseteq B \Big\}.
\end{array}
$$

\begin{claim} \label{bemerk}
We have $\gamma \leq N^2$, hence
there exists a nonempty open subset $W \subseteq B$ such that
$\LL \otimes \BB^N \otimes \KK^{k+1} \otimes \HH^{N^2}$ is weakly 
positive over $U \cap p^{-1}(W)$.
\end{claim}

\pproof \ref{bemerk}: $\HH^N \otimes \KK^{k+1}$ is ample with respect to
$U$.
So by \ref{amplunabh} there exists some $\beta > 0$ such that
$$\left(\LL \otimes \BB^N \otimes \KK^{k+1} \otimes
\HH^{\gamma}\right)^{(N-1) \cdot \beta} \otimes \left( \HH^N \otimes
\KK^{k+1} \right)^{\beta}, $$ and hence $$\left(\LL \otimes
\BB^{N-1} \otimes \KK^{k+1} \otimes \HH^{\frac{\gamma \cdot
(N-1)+N}{N}} \right)^{N \cdot \beta} \otimes \LL^{- \beta} $$ is
globally generated over $U \cap p^{-1}(W)$. We may assume that 
$\LL^{\beta}(-\beta \cdot D)$ is globally generated over $U$. 
Taking $$\MM = \LL \otimes \BB^{N-1}  \otimes \HH^{\frac{\gamma \cdot
(N-1)+N}{N}},$$ we obtain that $$(\MM \otimes \KK^{k+1})^{N
\cdot \beta} = \OO_Z(\Gamma + \beta \cdot D)$$ for some general
section $\Gamma$. Since $p$ is flat,
we may assume that $$p_{\ast} \left(\MM \otimes
\KK^{k+1} \otimes \oideal{Z/B}{\Gamma + \beta \cdot D}{\beta \cdot
N} \right)$$ is locally free over $W$, and hence, by \ref{versmult}, 
weakly positive over $W$. By \cite{ev84}, 4.3 we obtain that
$$p^{\ast} p_{\ast} \left(\MM \otimes \KK^{k+1} \otimes
\oideal{Z/B}{\Gamma + \beta \cdot D}{\beta \cdot N} \right)$$ is
weakly positive  over $p^{-1}(W)$, too.

The sheaves $\MM$ and $\KK$ satisfy the
conditions of \ref{surj}. Moreover, making $\beta$ larger if required, 
we may assume 
$$\oideal{Z/B}{\Gamma + \beta \cdot D}{\beta \cdot N} =
\oideal{Z/B} {D}{N}.$$ 
Thus we obtain that the composed map 
\begin{center}
\begin{tabular}{l}
$p^{\ast} p_{\ast} \left(\MM \otimes \KK^{k+1}
\otimes \oideal{Z/B} {\Gamma + \beta \cdot D}{\beta \cdot N}
\right) \longrightarrow$ \\ \\
$\MM \otimes \KK^{k+1} \otimes
\oideal{Z/B}{\Gamma + \beta \cdot D} {\beta \cdot N} \longrightarrow \MM
\otimes \KK^{k+1} \otimes \omega_{Z/B}$ 
\end{tabular}
\end{center}
is surjective over $U
\cap p^{-1}(W)$, since $e(D|_U) \leq N$. 
Using \cite{v95}, 2.16, we obtain that $\MM
\otimes \KK^{k+1} \otimes \omega_{Z/B}$ is weakly positive over 
$U \cap p^{-1}(W)$,
and so is
$$(\MM \otimes \KK^{k+1} \otimes \omega_{Z/B}) = \MM \otimes
\KK^{k+1} \otimes \bigoplus^r \omega_{Z/B}.$$ 
Again by \cite{v95}, 2.16 we obtain the weak positivity of
the quotient sheaf $$\MM \otimes \KK^{k+1} \otimes \BB =  \LL
\otimes \BB^N \otimes \KK^{k+1} \otimes \HH^{\frac{\gamma \cdot
(N-1)+N}{N}}$$ over $U \cap p^{-1}(W)$. But by definition of 
$\gamma$ we obtain
$$\frac{\gamma \cdot (N-1) + N}{N} > \gamma - N ,$$ and so
$\gamma \leq N^2$. \bewend

Now we have to get rid of the twist $\HH^{N^2}$, which we manage
to do by \ref{ueberh} and the following fact:

\begin{lemma}\label{vorueber}
In the situation above, let $\tau : B' \to B$ be a finite
morphism, ramified over some divisor $\Delta \subseteq B$, where
$D$ intersects the divisor $p^{-1}(\Delta)$ in codimension $\geq 2$. 
Let $Z' = B' \times_B Z$ denote the fibre product. 
Then, after replacing $B'$ by the complement of a closed subvariety 
of codimension $\geq 2$ if necessary, the assumptions made in
\ref{prop} 
hold true for the induced morphism $p' : Z' \to B'$ as well.
\end{lemma}

\proof By assumption we may assume that $D \cap p^{-1}(\Delta)$ is
even empty and that the induced morphism $\tau' : Z' \to Z$ 
is smooth over $Z -
p^{-1}(\Delta)$. Since moreover $p' : Z' \to B'$ is smooth over
$\tau^{-1}(B - \Delta)$, we obtain that $Z'$ is smooth and so in
particular $\omega_{Z'/B'} = {\tau'}^{\ast}\omega_{Z/B}$. The
compatibility of the relation $\varepsilon \geq 1$ with the
covering holds by \cite{f84}, 4.3. It remains to show that
the relation $e(D|_U) \leq N$ still holds. But taking $U' =
{\tau'}^{-1}(U)$ and $\Delta' = \tau^{\ast}\Delta$, we may assume, 
after making $W'$ smaller if necessary, that $U' \cap
{p'}^{-1}(\Delta')$ is empty and the relation then follows from
\ref{multetal}. \bewend

We now choose for every $\delta \in \N-\{0\}$ a finite morphism
$\tau_{\delta} : B_{\delta} \to B$, which satisfies the properties
of \ref{vorueber}. In addition we may, by \cite{v95}, 2.1, assume
that $\tau_{\delta}^{\ast}\tilde{\HH} = {\tilde{\HH'}}^{\delta}$
for some invertible sheaf $\tilde{\HH'}$ on $B_{\delta}$. This
induces, for every $\delta \in \N-\{0\}$, a morphism
$\tau_{\delta}' : Z_{\delta} \to Z$ such that
${\tau_{\delta}'}^{\ast}\HH = {\HH'}^{\delta}$ for some invertible
sheaf $\HH'$ on $Z_{\delta}$. By \ref{vorueber},
the bound $\gamma \leq N^2$ holds for $p_{\delta} :
Z_{\delta} \to B_{\delta}$ as well. By \ref{ueberh} this implies,
taking $\mu = N^2$, the weak positivity of  $\LL \otimes \BB^N
\otimes \KK^{k+1}$ over $U \cap p^{-1}(W)$, which proves \ref{prop}
for the flat case.

Now let us sketch how to reduce to the flat case if $p : Z \to B$ 
is not flat. Let $p_0 : Z_0 \to B_0$ be the flat locus. 
Let $H$ denote the Hilbert scheme parametrizing the flat subvarieties
of $Z$ (see \cite{v95}, p. 42, for example) and inducing, by its
universal property, a rational map $B_0 \to H$. Extending this  by 
\cite{h77}, II.7.17.3 to a
morphism $B' \to H$, we obtain a factorization of the inclusion $B_0
\to B$ over some birational map $\sigma : B' \to B$. By \cite{h77},
II.9.8.1, the morphism $B' \to H$ corresponds to a flat morphism $p_0'
: Z_0' \to B_0'$, where $Z_0'$ turns out to be a component of the
fibre product $Z' = Z \times_{B} B'$.

\begin{claim}\label{bem1}
Let $\delta : Z''_0 \to Z'_0$ be a desingularization such
that the preimage of the singular locus of $Z'_0$ is a
divisor and let ${\sigma_0''} : Z''_0 \to Z$ and $p''_0 : Z''_0 \to
B'$ denote the induced maps.  
Then, replacing the morphism $p : Z \to B$ in \ref{prop} by
$$p_0'' : Z''_0 - {p_0''}^{-1}\left(p_0''(B)\right) \to B' -
p_0''(B),$$ where $B$ is the maximal divisor in $Z''_0$ such that
$\codim (p_0''(B)) \geq 2$, the assumptions from
\ref{prop} hold true for the induced sheaves and divisors as well.
\end{claim}

To prove \ref{bem1}, the essential fact is the existence of a morphism 
$\omega_{Z_0''/Z} \to {p_0''}^{\ast}\omega_{B'/B}$ or, correspondingly, 
$\omega_{Z_0''/B'} \to {\sigma_0''}^{\ast}\omega_{Z/B}$. But this
follows from duality of finite morphisms (see \cite{h77},
III. Exercise 6.10) and from the fact that 
${p'}^{\ast}\omega_{B'/B} \simeq \omega_{Z'/Z}$ or, correspondingly, 
${\sigma'}^{\ast}\omega_{Z/B} \simeq \omega_{Z'/B'}$, which can be
proved using methods from \cite{h77}, III.6 and 7.

By \cite{h77}, III.10.2, we find
that $p''_0 : Z''_0 \to B'$, being an equidimensional morphism of
smooth varieties outside $B$, is flat outside $B$, which by
\ref{bem1} and \ref{prop} for the flat case implies the weak
positivity of ${\sigma_0''}^{\ast}\left(\LL \otimes
\BB^N \otimes \KK^{k+1}\right)|_{Z'' - B}$ over ${\sigma_0''}^{-1}(U) 
\cap
{p''_0}^{-1}(W')$ for an open subset $W' \subseteq B'$. We may
assume that ${\sigma_0''}^{-1}(U) \cap {p''_0}^{-1}(W')$ is of the
form ${\sigma_0''}^{-1}(U \cap p^{-1}(W))$ for an open subset $W
\subseteq B$. But then $\LL \otimes \BB^N \otimes 
\KK^{k+1}$ is weakly positive
over $U \cap p^{-1}(W)$ by \cite{ev84}, 4.3, which completes the proof of
\ref{prop}. \bewend


\begin{center}
\section{\sc A covering construction}
\end{center}

Following the notations from the introduction, 
we shall construct a covering which simplifies the index
conditions defining the ideal sheaf $\II$.

\begin{construction}\label{konstrue}
\rm
Let $\nu \in \{1 , \dots , n\}$.
For every point $\xi_{\mu , \nu} \in p_{\nu}(S)$ we choose exactly
$m_{\nu}$ smooth  hyperplane sections $H^{\nu}_{\mu,1} , \dots ,
H^{\nu}_{\mu,m_{\nu}}$ of $\LL_{\nu}$ in general position through
$\xi_{\mu , \nu}$. 
(If $|p_{\nu}(S)| = 1$ we just add another arbitrary point.)
So $\xi_{\mu , \nu} \in \cap_{k = 1}^{m_{\nu}}
H^{\nu}_{\mu,k}$ is an isolated point, where the hyperplane
sections $H^{\nu}_{\mu,k}$ intersect transversally. We denote by
$$\Delta_{\nu} = \sum H^{\nu}_{\mu,k}$$ the corresponding normal crossing
divisor on $X_{\nu}$, where the sum runs over $\mu =
1 , \dots , M_{\nu}'$ and $k = 1 , \dots , m_{\nu}$. So we have
$$\OO_{X_{\nu}}(\Delta_{\nu}) = \LL_{\nu}^{M_{\nu}' \cdot
m_{\nu}}.$$ Now by \cite{k81}, Theorem 17, there exists a
smooth projective variety $Y_{\nu}$ and a finite morphism
$\sigma_{\nu} : Y_{\nu} \to X_{\nu}$ such that
$\sigma_{\nu}^{\ast}\Delta_{\nu}$ is a normal crossing divisor on
$Y_{\nu}$ and $\sigma_{\nu}$ ramifies exactly over $\Delta_{\nu}$,
where, taking $$N = \min \Big\{n \in \N-\{0\} \, ; \,
\frac{n}{d_{\nu}} \cdot t_{\mu} \in \N-\{0\} \, \mbox{for} \, \mu
=1, \dots , M  \, \mbox{and} \, \nu =1, \dots , n \Big\},$$ for
$\mu = 1 , \dots , M$ and $k = 1 , \dots , m_{\nu}$ we have:
$$\sigma_{\nu}^{\ast} H^{\nu}_{\mu,k} = \frac{N}{d_{\nu}} \cdot
\left(\sigma_{\nu}^{\ast} H^{\nu}_{\mu,k}\right)_{red}.$$ The
coverings $\sigma_{\nu} : Y_{\nu} \to X_{\nu}$ constructed in this
way for $\nu = 1 , \dots , n$ induce a covering $$\sigma : Y = Y_1
\times \dots \times Y_n \longrightarrow X,$$
which is \'etale over 
$X_0 = (X_1-\Delta_1) \times \dots \times (X_n - \Delta_n)$.
\end{construction}

\begin{lemma}\label{koro1}
In the situation of \ref{konstrue} we consider the sheaf $$\MM =
\bigcap_{\mu = 1}^M \bigcap_{\eta \in \sigma^{-1}(\xi_{\mu})}
m_{\eta}^{N \cdot t_{\mu}}, $$ which is full. Let
${\tau'} : Z \to Y$ be a birational morphism such
that $Z$ is smooth and $\MM' = {\tau'}^{-1}\MM \cdot \OO_{Z}$ is
invertible on $Z$ and 
let us fix a blowing up $\tau : X' \to X$ such that $\II' =
\tau^{-1}\II \cdot \OO_{X'}$ is invertible. Let us assume
moreover that there exists a morphism $\sigma' :
Z \to X'$ making the diagram
\begin{eqnarray*}
Z & \stackrel{{\tau'}}{\longrightarrow} & Y \\
\sigma' \downarrow & & \downarrow \sigma\\
X' & \stackrel{\tau}{\longrightarrow} & X .
\end{eqnarray*}
commutative. Let us denote by $\varrho : Z \to X$ the induced 
morphism. Then for all $l \in \N-\{0\}$ the
trace map induces a surjective map
$\varrho_{\ast} {\MM'}^l \to \II^l$, and moreover,
the ideal sheaf $\II^l$ is full.
\end{lemma}

\proof Keeping notations simple, we restrict ourselves to the case
$l = 1$. For the general case one has to consider the $l$-th powers.
We shall first show that the image of
$\sigma_{\ast}\MM$ under the trace map $\sigma_{\ast}\OO_Y \to
\OO_X$ is contained in $\II$. This statement is local, and we may,
after fixing a point $\xi_{\mu} \in S$ and some point $\eta \in Y$
mapping to $\xi_{\mu}$
and keeping all the
notations, replace $X$ and $Y$ by
the corresponding local rings.

Let $h^{\nu}_{\mu,k}$ and $g_{\mu,k}^{\nu}$ denote the local equations 
of $H^{\nu}_{\mu,k}$ in
the point $\xi_{\mu , \nu}$ and of $(\sigma_{\nu}^{\ast}
H^{\nu}_{\mu,k})_{red}$ in a point $\eta$ mapping to $\xi_{\mu,\nu}$,
respectively. 
Thus for $\mu = 1 , \dots , M$, the ideal sheaves 
$\II_{\xi_{\mu}, t_{\mu}}$ and $\sigma^{\ast}\II_{\xi_{\mu} ,
  t_{\mu}}$ are 
generated by expressions of the form $$\prod_{k = 1}^{m_1}
(h^1_{\mu,k})^{\alpha_1} \cdots \prod_{k = 1}^{m_n} (h^n_{\mu ,
k})^{\alpha_n}$$ where $\frac{\alpha_1}{d_1} + \dots +
\frac{\alpha_n}{d_n} \geq t_{\mu}$ and $$\prod_{k = 1}^{m_1} (g_{\mu ,
k}^{1})^{\frac{\alpha_1 \cdot N} {d_1}} \cdots \prod_{k = 1}^{m_n}
(g_{\mu , k}^{n})^{\frac{\alpha_n \cdot N} {d_n}} $$ where
$\frac{\alpha_1 \cdot N}{d_1} + \dots + \frac{\alpha_n \cdot
N}{d_n} \geq N \cdot t_{\mu},$ respectively. 
It remains to show that the image of $\sigma_{\ast}m_{\eta}^{N
\cdot t_{\mu}}$ under the trace map is contained in
$\II_{\xi_{\mu} , t_{\mu}}$. The coverings $\sigma_{\nu} : Y_{\nu}
\to X_{\nu}$ are determined by the ramifications $\beta_{\nu} =
\frac{N}{d_{\nu}} \in \N-\{0\}$ in $\Delta_{\nu}$, and we may now
assume that $\sigma : Y \to X$ is a Galois covering with Galois
group $G = \Z/\beta_1 \times \dots \times \Z/\beta_n$. The trace
map then is nothing but the sum over all conjugates of $G$. Now
$m_{\eta}^{N \cdot t_{\mu}}$ is invariant under $G$, so the image
of $\sigma_{\ast}m_{\eta}^{N \cdot t_{\mu}}$ under the trace map
in $\OO_X$ is generated by the images of the $G$-invariant
elements $$\sum_{\underline{l} = (l_1 , \dots , l_n)}
f_{\underline{l}} \cdot \left(\prod_{k=1}^{m_1}
g_{\mu,k}^1\right)^{l_1} \cdots \left(\prod_{k=1}^{m_n}
g_{\mu,k}^n\right)^{l_n},$$ where $f_{\underline{l}}$ are units for
$\underline{l} \in \N^n$. But the $G$-invariance of such expressions
just means the $G$-invariance of the single summands, and
for such a summand, given by $\underline{l} = (l_1 , \dots , l_n)$, 
this just means that
$l_{\nu} = \beta_{\nu} \cdot s_{\nu}$ for $\nu = 1 , \dots , n$
and certain $s_{\nu} \in \N-\{0\}$. So the image of
$\sigma_{\ast}m_{\eta}^{N \cdot t_{\mu}}$ under the trace map in
$\OO_X$ is generated by expressions of the form
$$\left(\prod_{k=1}^{m_1} h_{\mu,k}^1\right)^{s_1} \cdots
\left(\prod_{k=1}^{m_n} h_{\mu,k}^n\right)^{s_n},$$ where $\sum_{\nu
= 1}^N l_{\nu} = \sum_{\nu = 1}^N \beta_{\nu} \cdot s_{\nu} \geq N
\cdot t_{\mu}$. This means that $\sum_{\nu = 1}^N
\frac{s_{\nu}}{d_{\nu}} \geq t_{\mu},$ so the image of
$\sigma_{\ast}m_{\eta}^{N \cdot t_{\mu}}$ under the trace map is
contained in $\II_{\xi_{\mu} , t_{\mu}}$.

So the trace map induces a surjective map
$\sigma_{\ast}\MM \to \II$. Since $\MM$ is full, we have
$\tau'_{\ast}\MM' = \MM$, and we obtain a surjective map
$\varrho_{\ast}\MM' \to \II$.

Finally, we observe that the inclusion
$\tau_{\ast}{\II'} \to \sigma_{\ast}{\tau'}_{\ast}{\MM'} =
\sigma_{\ast}\MM$ splits, and we obtain a surjective map
$\sigma_{\ast}\MM \to \tau_{\ast}{\II'}$, which factors over
$\II$, as we have just seen. 
So $\II \to \tau_{\ast}{\II'}$ is surjective and $\II$ is 
full. \bewend


\begin{center}
\section{\sc The Proof of \ref{dysonv}}
\end{center}

We shall now give the proof for \ref{dysonv}.
First we need some more notations.

\begin{definition}
\rm An open subset $U \subseteq X$ is called {\it open subset of
type $k$}, if for $\nu = 1 , \dots , k$ there exist open subsets
$U_{\nu} \subseteq X_{\nu}$ and an open subset $W \subseteq
X_{k+1} \times \dots \times X_n$ such that $U$ is of the form $U =
U_1 \times \dots \times U_k \times W$. A subset of type $n-1$ is
called a {\it product open set}.
\end{definition}

As one immediately sees, \ref{dysonv} follows from inductively
applying the following

\begin{theorem}\label{mainlemma}
If in the situation of \ref{dysonv} the sheaf
$\underline{\LL}^{\underline{d}} \otimes \II$ is weakly  positive
over some open set $U^{k-1}$ of type $k-1$, then the sheaf
$\underline{\LL}^{\underline{d'}} \otimes \II$ is weakly positive
over some open set $U^k$ of type $k$, where $\underline{d'} = (d_1'
, \dots , d_n')$ with
$$
d_{\nu}' = \left\{ \begin{array}{ll}
d_{\nu} + M_{\nu} \cdot d_{k+1} \cdot e_{k+1} & (\nu \leq k) \\
d_{\nu} & (\nu \geq k+1)
\end{array}
\right.
$$
and $M_{\nu} = M_{\nu}' \cdot m_{\nu} + \gamma_{\nu}$ for $\nu = 1
, \dots , n$.
\end{theorem}

So we shall prove this. To this end let $\alpha > 0$. 
By \ref{ivoll3} there
exists some $\beta > 0$ such that
$\underline{\LL}^{\underline{\delta}} \otimes \JJ$ is globally
generated over $U^{k-1} =
U_1 \times \dots \times U_{k-1} \times W$,
where $\underline{\delta} = (\delta_1 , \dots , \delta_n)$ with
$\delta_{\nu} = \alpha \beta d_{\nu} + \beta$ for $\nu = 1 , \dots
, n$ and $\JJ = \II^{\alpha \beta}$ and where $U_{\nu} \subseteq
X_{\nu}$ for $\nu = 1 , \dots , k-1$ and $W \subseteq U_k \times
\dots \times U_n$ are open subsets. We notice that the ordering of
the $d_{\nu} \cdot e_{\nu}$ is preserved; we have $\delta_1 \cdot
e_1 \geq \dots \geq \delta_n \cdot e_n$. We choose a general
section $\Gamma$ of $\underline{\LL}^{\underline{\delta}} \otimes
\JJ$ and observe:

\begin{claim}\label{absch}
There exists a product open subset $U' \subseteq X$ such that
$$e(\Gamma|_{U'}) \leq \delta_{k+1} \cdot e_{k+1} + 1.$$ Moreover,
by \ref{konstrue} we may assume that
$p|_{p^{-1}(U')}$ is \'etale, and so by \ref{multetal}
we have $e(p^{\ast}\Gamma|_{p^{-1}(U')}) \leq
\delta_{k+1} \cdot e_{k+1} + 1$.
\end{claim}

\proof We consider the projection $p_{1 , \dots , k} : X \longrightarrow
X_1 \times \dots \times X_k$ onto the first $k$ factors of the
product, whose fibres are isomorphic to $X_{k+1} \times \dots
\times X_n$. One has $p_{1 , \dots , k}(U^{k-1}) = U_1 \times
\dots \times U_{k-1} \times U_k$ for an open subset $U_k \subseteq
X_k$. Let $x \in p_{1 , \dots , k}(U^{k-1})$, in other words, the
fibre $p_{1 , \dots , k}^{-1}(x)$ intersects the open subset
$U^{k-1}$. Then we obtain $$e(\Gamma|_{p_{1 , \dots , k}^{-1}(x)})
\leq e(\underline{\LL}^{\underline {\delta}}|_{p_{1 , \dots ,
k}^{-1}(x)}) = e(\LL_{k+1}^{\delta_{k+1}} \boxtimes \dots \boxtimes
\LL_n^{\delta_n}).$$ 
By \cite{v95}, 5.11 we have $$e(\LL_{\nu}^{\delta_{\nu}})
\leq \delta_{\nu} \cdot c_1(\LL_{\nu})^{m_{\nu}} + 1 =
\delta_{\nu} \cdot e_{\nu} + 1$$ for $\nu = k+1 , \dots , n$.
Now by \ref{products} and since the $\delta_{\nu} \cdot e_{\nu}$ are
ordered, we obtain
$$e(\Gamma|_{p_{1 , \dots , k}^{-1}(x)}) \leq \delta_{k+1} \cdot
e_{k+1} +1,$$ which, by \cite{v95}, 5.14, yields
$e(\Gamma|_{U'}) \leq \delta_{k+1} \cdot e_{k+1} + 1$, where
$U'$ is a neighbourhood of the fibre $p_{1 , \dots , k}^{-1}(x)$.
We can say more precisely that $U' = p_{1 , \dots , k}(U^{k-1})
\times X_{k+1} \times \dots \times X_n \cap X_0,$ and thus is a
product open subset. \bewend
 
In order to apply \ref{prop} to our situation, we shall now check
whether the assumptions are fulfilled.
For $\nu = 1 , \dots , n$ let $\sigma_{\nu}
: Y_{\nu} \to X_{\nu}$ be the coverings constructed in
\ref{konstrue}, which ramify exactly in $\Delta_{\nu}$, and let
$\Delta_{\nu}' = p_{\nu}^{-1}(\Delta_{\nu})$. Let 
$\Delta = \sum_{\nu = 1}^{k} \Delta_{\nu}'$ and let us write 
$\sigma_{1 , \dots , k} = \sigma_1 \times \dots \times \sigma_k$
and $\sigma_{k+1 , \dots , n} = \sigma_{k+1} \times \dots \times
\sigma_n$.

\begin{claim}\label{omega}
In the situation described above, we have
$$\omega_{Y/(Y_{k+1} \times \dots \times
Y_n)} = \sigma^{\ast} \omega_{X/(X_{k+1} \times \dots \times X_n)}
\otimes \OO_Y\left(-(\sigma^{\ast}\Delta)_{red} +
\sigma^{\ast}\Delta\right).$$
\end{claim}

\pproof \ref{omega}: 
We consider the fibre product
\begin{eqnarray*}
X_1 \times \dots \times X_k
\times Y_{k+1} \times \dots \times Y_n & \stackrel{\sigma'_{k+1,\dots,n}}
{\longrightarrow} & X \\
\downarrow & & \downarrow \\
Y_{k+1} \times \dots \times Y_n &
\stackrel{\sigma_{k+1,\dots,n}}{\longrightarrow}
& X_{k+1} \times \dots \times X_n.
\end{eqnarray*}
and notice that
$\sigma'_{k+1,\dots,n}$ ramifies exactly in $\sum_{\nu = k+1}^n
\Delta'_{\nu}$. Compatibility of relatively canonical sheaves with
base change yields $$\omega_{(X_1 \times \dots \times X_k \times
Y_{k+1} \times \dots \times Y_n)/ (Y_{k+1} \times \dots \times
Y_n) } = {\sigma'_{k+1,\dots,n}}^{\ast}\omega_{X/(X_{k+1} \times
\dots \times X_n)},$$ and this holds true after adding divisors,
which are not contained in the ramification locus, hence we have 
$$\omega_{(X_1
\times \dots \times X_k \times Y_{k+1} \times \dots \times Y_n)/
(Y_{k+1} \times \dots \times Y_n) } \left( \log
{\sigma'_{k+1,\dots,n}}^{\ast} \left(\sum _{\nu = 1}^k
\Delta'_{\nu}\right)\right)$$ 
$$ = {\sigma'_{k+1,\dots,n}}^{\ast}\left(\omega_{X/(X_{k+1} \times
\dots \times X_n)} \left(\log\left(\sum_{\nu=1}^k
\Delta'_{\nu}\right) \right)\right).$$ 
Considering on the other hand the fibre product
\begin{eqnarray*}
Y & \stackrel{\sigma'_{1,\dots,k}}{\longrightarrow} & X_1 \times \dots
\times X_k
\times Y_{k+1} \times \dots \times Y_n  \\
\downarrow & & \downarrow  \\
Y_1 \times \dots \times Y_k &
\stackrel{\sigma_{1,\dots,k}}{\longrightarrow} &
X_1 \times \dots \times X_k,  \\
\end{eqnarray*}
then $\sigma'_{1,\dots,k}$ ramifies exactly in $\sum_{\nu = 1}^k
{\sigma'_{k+1,\dots,n}}^{\ast} \Delta'_{\nu}$, and the Riemann
Hurwitz Formula implies $$\omega_{Y/(Y_{k+1} \times \dots \times
Y_n)}\left(\log {\sigma'_{1,\dots,k}}^{\ast} \left(\sum_{\nu =
1}^k {\sigma'_{k+1,\dots,n}}^{\ast} \Delta'_{\nu}\right)\right)$$
$$= {\sigma'_{1,\dots,k}}^{\ast}\left(\omega_{(X_1 \times \dots
\times X_k \times Y_{k+1} \times \dots \times Y_n)/(Y_{k+1} \times
\dots \times Y_n) } \left(\log\left(\sum_{\nu = 1}^k
{\sigma'_{k+1,\dots,n}}^{\ast}\Delta'_{\nu}\right)\right)\right).$$
$$={\sigma'_{1,\dots,k}}^{\ast}\left({\sigma'_{k+1,\dots,n}}^{\ast}
\left(\omega_{X/(X_{k+1}
\times \dots \times X_n)} \left(\log\left(\sum_{\nu=1}^k
\Delta'_{\nu}\right) \right)\right)\right).$$ Thus we get
$$\omega_{Y/(Y_{k+1} \times \dots \times Y_n)}\left(\log
\sigma^{\ast}\left(\sum_{\nu = 1}^k \Delta'_{\nu}\right)\right)$$
$$= \sigma^{\ast}\left(\omega_{X/(X_{k+1} \times \dots \times
X_n)} \left(\log\left(\sum_{\nu=1}^k \Delta'_{\nu}\right)
\right)\right)$$ and so $$\omega_{Y/(Y_{k+1} \times \dots \times
Y_n)} \otimes \OO_Y\left(\left( \sigma^{\ast}\left(\sum_{\nu =
1}^k \Delta'_{\nu}\right)\right)_{red}\right)$$
$$=\sigma^{\ast}\left(\omega_{X/(X_{k+1} \times \dots \times X_n)}
\otimes \OO_X\left(\left(\sum_{\nu = 1}^k
\Delta'_{\nu}\right)_{red}\right)\right).$$  
This proves \ref{omega}. \bewend

Now $\tau' : Z \to Y$ is the blowing
up of the ideal sheaf $\MM$, given by powers of maximal ideals, so
we have $$\omega_{{Z}/(Y_{k+1} \times \dots \times Y_n)} =
{\tau'}^{\ast}(\omega_{Y/(Y_{k+1} \times \dots \times Y_n)})
\otimes \OO_{{Z}}(E)$$ $$ = p^{\ast}\omega_{X/(X_{k+1}
\times \dots \times X_n)} \otimes p^{\ast}\OO_X(\Delta)
\otimes {\tau'}^{\ast}\OO_Y \left(-(\sigma^{\ast}\Delta)_{red}\right)
\otimes \OO_{{Z}}(E),$$ where $E$ is the exceptional divisor of
${\tau'}$. By choice of $\gamma_{\nu}$ there exists a surjective map
$$\bigoplus^r \omega_{X/(X_{k+1} \times \dots
\times X_n)} \longrightarrow \underline{\LL}^{(\gamma_1 , \dots ,
\gamma_k , 0 , \dots , 0)},$$ so we obtain a surjective map
\begin{eqnarray*}
\bigoplus^r \omega_{{Z}/(Y_{k+1} \times \dots \times Y_n)} & \to &
\varrho^{\ast} \underline{\LL}^{(\gamma_1 , \dots , \gamma_k , 0 , \dots ,
0)}
\otimes \varrho^{\ast} \underline{\LL}^{(M_1' \cdot m_1 , \dots , M_k'
\cdot m_k ,
0 , \dots , 0)} \\
& & \otimes \, {\tau'}^{\ast}\OO_Y
\left(-(\sigma^{\ast}\Delta)_{red}\right) \otimes
\OO_{{Z}}(E) \\
& = & \varrho^{\ast} \underline{\LL}^{(M_1 , \dots , M_k , 0 , \dots , 0)}
\\
& & \otimes \, {\tau'}^{\ast}\OO_Y
\left(-(\sigma^{\ast}\Delta)_{red}\right) \otimes
\OO_{{Z}}(E),
\end{eqnarray*}
where the last sheaf coincides with $\varrho^{\ast}
\underline{\LL}^{(M_1 , \dots , M_k , 0 , \dots , 0)}$ over
$\varrho^{-1}(U')$. Hence,
there exists a map $$\bigoplus^r \omega_{{Z}/(Y_{k+1} \times
\dots \times Y_n)} \longrightarrow
\varrho^{\ast}\underline{\LL}^{(M_1 , \dots , M_k , 0 , \dots ,
0)},$$ surjective over $\varrho^{-1}(U')$.

Now for $B = Y_{k+1} \times \dots \times Y_n$ and taking
$\LL = \varrho^{\ast}\underline{\LL}^{\underline{\delta}}$, 
$\BB = \varrho^{\ast}\underline{\LL}^{(M_1 , \dots , M_k , 0 , \dots ,
  0)}$ and 
$\KK = \varrho^{\ast}\underline{\LL}^{(1 , \dots , 1)}$,
we may apply \ref{prop} to the induced morphism 
$p: Z \to B$ and obtain that 
$\varrho^{\ast}\underline{\LL}^{\underline{\delta'}}
\otimes {\MM'}^{\alpha \beta}$ is weakly positive over 
$$\varrho^{-1}(U') \cap
\varrho^{-1}(W) = \varrho^{-1}(U') \cap {\tau'}^{-1} (Y_1 \times \dots
\times Y_k \times W),$$

where  $W \subseteq
Y_{k+1} \times \dots \times Y_n$ is a non-empty subset and 
$\underline{\delta'} = (\delta'_1 , \dots , \delta'_n)$ with
$$
\delta_{\nu}' = \left\{ \begin{array}{ll}
\delta_{\nu} + M_{\nu} \cdot (\delta_{k+1} \cdot e_{k+1} + 1) + k + 1 &
(\nu \leq k) \\
\delta_{\nu} + k + 1 & (\nu \geq k+1).
\end{array}
\right.
$$

Now we may assume that $\varrho^{-1}(U') \cap
{\tau'}^{-1} (Y_1 \times \dots  \times Y_k \times W)$ is of the
form $\varrho^{-1}(U^k)$, where $U^k \subseteq X$ is an open
subset of type $k$. By \ref{koro1}, for every $l \in \N -
\{0\}$ there exists a surjective map
$$\varrho_{\ast}\left(\varrho^{\ast}\underline{\LL}^{\underline{\delta'}}
\otimes {\MM'}^{\alpha \beta}\right)^l =
(\underline{\LL}^{\underline{\delta'}})^l \otimes
\varrho_{\ast}{\MM'}^{\alpha \beta l} \to
(\underline{\LL}^{\underline{\delta'}})^l \otimes \JJ^l, $$ and 
$\JJ^l$ is full for every $l \in \N - \{0\}$. So
taking $$\GG = \varrho^{\ast}\underline{\LL}^{\underline{\delta'}}
\otimes {\MM'}^{\alpha \beta} \,\, {\mbox {and}} \,\, \LL =
\underline{\LL}^{\underline{\delta'}},$$ 
and since $\varrho|_{\varrho^{-1}(U^k)}$ is finite, 
the assumptions of
\ref{posrunter} are fulfilled, and we obtain the weak positivity
of $\underline{\LL}^{\underline{\delta'}} \otimes \JJ$ over $U^k$.
Now we have
\begin{eqnarray*}
\underline{\LL}^{\underline{\delta'}} \otimes \JJ & = &
\left(\underline{\LL}^{(d_1
+ M_1 \cdot d_{k+1} \cdot e_{k+1} , \dots , d_k + M_k \cdot d_{k+1}
\cdot e_{k+1} , d_{k+1}, \dots , d_n)}\right)^{\alpha \beta} \otimes
\II^{\alpha \beta} \\
& & \otimes \, \left(\underline{\LL}^{(M_1 \cdot e_{k+1} + 1 , \dots , M_k
\cdot e_{k+1}+ 1 , 1 , \dots , 1)}\right)^{\beta} \\
& & \otimes \, \underline{\LL}^{(M_1 + k + 1 , \dots , M_k + k+1 , k+1 ,
\dots , k+1)}.
\end{eqnarray*}
Hence by \ref{ivoll3} we obtain the weak positivity of
$$\underline{\LL}^{(d_1 + M_1\cdot d_{k+1} \cdot e_{k+1} , \dots ,
d_k + M_k \cdot d_{k+1} \cdot e_{k+1} , d_{k+1}, \dots , d_n)}
\otimes \II = \underline{\LL}^{\underline{d'}} \otimes \II$$ over
the open set $U^k \subseteq X$ of type $k$, which completes the proof
of \ref{mainlemma}. \bewend


\begin{center}
\section{\sc The Curve Case}
\end{center}

As explained in the introduction, the case of curves is, in some
sense, the most 
``natural'' one. Let us hence evaluate \ref{dysonv} for the special 
situation of a product 
$C = C_1 \times \dots \times C_n$ of smooth projective curves of genera
$g_1 , \dots , g_n$ (which,
indeed, was the motivation for this paper). Using
notations from above and taking $\LL_{\nu} = \OO_{C_{\nu}}(1)$, hence 
$e_{\nu} = 1$ and $\gamma_{\nu} = 2 g_{\nu} - 2$ for $\nu = 1 ,
\dots , n$, we obtain $\underline{\LL}^{\underline{d}} =
\OO_C(\underline{d})$. 
Near a point $\xi \in C$, sections of this sheaf have an expansion
$$s = \sum_{\alpha_{\nu} \leq \delta_{\nu}} \lambda_{\alpha_1 , \dots ,
\alpha_n} \cdot
z_1^{\alpha_1} \cdots z_n^{\alpha_n},$$
for local parameters $z_1 , \dots , z_n$ in $\xi$.
Moreover, we can define in the usual way
the index $\text{ind}_{\xi}(s)$ of such a section $s$ in the point $\xi$
by
$$\text{ind}_{\xi}(s) = \min \left\{ \sum_{\nu = 1}^n
\frac{\alpha_{\nu}}{d_{\nu}} \,
; \, \lambda_{\alpha_1 , \dots , \alpha_n} \neq 0 \right\}.$$
The ideal sheaves
$\II_{\xi_{\mu} , t_{\mu}}$ turn out to be generated 
by global sections $s$ satisfying 
$\text{ind}_{\xi_{\mu}}(s) \geq t_{\mu}$.
In this situation we obtain from \ref{dysonv} the following result:

\begin{corollary}[Dyson's Lemma for curves]\label{dysonc2}
If $\OO_C(\underline{d}) \otimes \II$ is effective, then
$\OO_C(\underline{d'}) \otimes \II$ is weakly positive over a
product open set in $C$, where $\underline{d'} = (d_1' , \dots ,
d_n')$ with $$d'_{\nu} = d_{\nu} + M_{\nu} \cdot \sum_{j= \nu +
1}^n d_j$$ and $M_{\nu} = 2 g_{\nu} - 2 + M_{\nu}'$ for $\nu = 1 ,
\dots , n$.
\end{corollary}

\begin{remark}
One can do slightly better here. Making a very mild assumption on the
position of the points and modifying the ideal sheaves a little, one
can replace the notion of weak positivity over some product open set by 
the notion of numerical effectivity, hence by a global positivity 
statement. We shall, however, not prove this here, since it
can be achieved by exactly the same arguments used in \cite{ev84}, 5.

\end{remark}

\vspace{2cm} 

\scriptsize

Markus Wessler, Universit\"at Kassel, Fachbereich 17, 
Mathematik und Informatik

Heinrich-Plett-Str. 40, D-34132 Kassel

e-Mail: {\tt wessler@mathematik.uni-kassel.de}

\end{document}